\documentclass[amsfonts]{article}

\newtheorem{theorem}{Theorem}

\newtheorem{corollary}[theorem]{Corollary}

\newenvironment{proof}[1][Proof]{\textbf{#1.} }{\ \rule{0.5em}{0.5em}}

\newcommand\ccr{{\mathcal {R}}}

\newcommand\ccl{{\mathcal {L}}}

\title{A characterization of Morita equivalence pairs of quantales}

\author{Jan Paseka\footnote {The paper was prepared under the support of the
{\em Grant Agency of the Czech Republic} (GA \v CR 201/00/1466).}}

\begin{document}
\maketitle

\begin{abstract}
We characterize the pairs of
sup-lattices which occur as pairs of Morita equivalence
bimodules between
quantales
in terms of the mutual relation between the sup-lattices.
\end{abstract}


{\bf MS classification:} 46M15, 46L05, 18D20, 06F07.


\section{Introduction and preliminaries}

In this paper we give a characterization of those pairs $(X,Y)$ of
sup-lattices, which occur as pairs of Morita equivalence bimodules between
m-regular quantales,
in terms of the mutual relation between the sup-lattices $X$ and $Y$.
The sup-lattice characterization is based on a surprising  similarity with
the C$^{*}$-algebra theoretical characterization of Morita pairs
given in \cite{todorov}.

Our main motivation comes from both Algebra and Analysis, where the
concept of Morita equivalence has been applied to many different
categories to explore the relationship between objects and their
representation theory i.e. the theory of modules.

The idea of Morita equivalence was first made precise by Morita \cite{morita}
in the context of the category of unital rings: two unital rings are called
{\em Morita equivalent}  if their categories of right (left) modules are
equivalent. Morita equivalent rings always come with a pair of
corresponding bimodules of a certain type in such a way that the
functors implementing the equivalence of the categories of modules are
actually equivalent to tensoring with these modules. Morita
equivalent rings share many ring theoretical properties - the
Morita invariants.

The notion of Morita equivalence has been adapted to many other
algebraic and analytic contexts such as nonunital rings \cite{anh-marki},
rings with involution \cite{ara}, C$^{*}$-algebras \cite{blecher},
\cite{rieffel}, monoids \cite{banaschewski}, unital quantales \cite{borceux}
and nonunital involutive quantales \cite{paseka}.

The paper is organized as follows. Section 1 introduces the basic
notions concerning the subject. Section 2 contains the basic
theorem for quantales, section 3 its counterpart for involutive quantales.
In section 4,  we show that our characterization results from section 1 and
section 2  translate to multiplicative (involutive) semilattices and, in
fact, to (involutive) $\kappa$-quantales.

\vskip0.1cm


A {\em quantale\/} is a sup-lattice $A$
with an associative
binary multiplication satisfying
$$
x\cdot\bigvee\limits_{i\in I}
x_i=\bigvee\limits_{i\in I}x\cdot
x_i\ \ \hbox{and}\ \ (\bigvee\limits_{i\in I}x_i)\cdot
x=\bigvee\limits_{i\in I}x_i\cdot x
$$
for all $x,\,x_i\in A,\,i\in I$ ($I$ is a set).
$1$ denotes the greatest element
of $A$, $0$ is the smallest element of $A$. A quantale
$A$ is said to be {\em unital} if there is an element
$e\in A$ such that $e\cdot a= a = a\cdot e$ for all
$a\in A$.


An {\em involution}
on a sup-lattice $S$ is a unary operation  such that
$$
\begin{array}{r c l}
a^{**}&=&a,\\
(\bigvee a_{i})^{*}&=&\bigvee a_{i}^{*}
\end{array}
$$
\noindent for all $a, a_{i}\in S$.
An {\em involution}
on a quantale $A$ is an involution on the sup-lattice
$A$  such that
$$
\begin{array}{r c l}
(a\cdot b)^{*}&=&b^{*}\cdot a^{*},\\
\end{array}
$$
\noindent for all $a,  b\in A$.
A sup-lattice (quantale) with the involution is said to be
{\em involutive}.


By a {\em morphism of} ({\em involutive}) {\em quantales}
will be meant a
$\bigvee$- ($^{*}$-) and $\cdot$\,{}-preserving mapping
$f:A\to A'.$ If a morphism preserves the unital element
we say that it is {\em unital}.


Let $A$ be a quantale.
A {\em right module over} $A$ (shortly a right $A$-module) is
a sup-lattice $M$, together with a {\em module action}
$\hbox{\rm{\_}}{\cdot}\hbox{\rm{\_}}:M\times A\to M$ satisfying

$$
\begin{array}{c l c l c}
\phantom{xxxxxxxxxx}&%
m\cdot(a\cdot b)&=&(m\cdot a)\cdot b&
\phantom{xxxxxxxxxx}\hbox{\rm (M1)}\\
\phantom{xxxxxxxxx}&%
(\bigvee X)\cdot a&=&\bigvee \{x\cdot a: x\in X\}&
\phantom{xxxxxxxxxx}\hbox{\rm (M2)}\\
\phantom{xxxxxxxx}&
m\cdot \bigvee S&=&\bigvee \{m\cdot s: s\in S\}&
\phantom{xxxxxxxxxx}\hbox{\rm (M3)}\\
\end{array}
$$

\noindent for all $a, b\in A$, $m\in M$, $S\subseteq A$,
$X\subseteq M$.

For a right $A$-module $X$
the submodule $\hbox{\rm ess}(X)=X\cdot A$
generated by the elements
$x\cdot a$ is called the {\em essential part} of $X$.
If $\hbox{\rm ess}(X)=X$ we say that $X$ is {\em essential}.

 We shall say that
$A$ is {\em right separating} for the $A$-module $M$
and
that $M$ is ({\em right}) {\em separated} by $A$ if
$m\cdot (-)=n\cdot (-)$ implies $m=n$.
We say that $M$ is {\em m-regular} if it is both separated by $A$
and essential.

All definitions and propositions stated for right $A$-modules are valid
in a dualized form for left $A$-modules. An ({\em m-regular}) {\em $A, B$-bimodule} is
a sup-lattice $M$, together with a left ({m-regular})  $A$-module action
and a right ({m-regular})
$B$-module action such that the left action with elements in $A$ and the
right action with elements in $B$ commute.

An (involutive) quantale $A$ is called
{\em m-regular} if it is {m-regular} as an $A, A$-bimodule.
 Then, evidently $1\cdot 1=1$ in $A$, $(-)\cdot a=(-)\cdot b$ implies $a=b$
and $a\cdot (-)=b\cdot (-)$ implies $a=b$.

For facts concerning quantales and quantale modules  in general we refer
to \cite{rosenthal1}.


To make our paper elementary and self-contained we will give
an explicit and elementary definition of Morita
equivalence of quantales (that is of course equivalent to
the usual one as it was recently shown by the author \cite{pasekam}).

Two m-regular quantales $A$ and $B$ are said to be {\em Morita equivalent} if
there exist sup-lattices $X$ and $Y$, such that

1)$X$ is an m-regular $A,B$-bimodule, $Y$ is an m-regular $B,A$-bimodule;

2) there are bimodule maps
$( -, -) : X\times Y\longrightarrow A$
and  $[ -, -] : Y\times X\longrightarrow B$ called pairings
such that $(x\cdot b,y) = (x, b\cdot y)$,
$[y\cdot a,x] = [y, a\cdot x]$ (that is, these maps are
{\em balanced}), $(x_1,y)\cdot x_2 = x_1\cdot [y,x_2]$,
$[y_1,x]\cdot y_2 = y_1\cdot (x,y_2)$ for each
$x,x_1,x_2\in X$, $y,y_1,y_2\in Y$, $a\in A$, $b\in B$ and

3) the sup-preserving maps on the  sup-lattice tensor
products   $X\otimes Y, Y\otimes X$
induced by the pairings are surjective.


The 6-tuple $(A,B,X,Y,( -, -),[ -, -])$ is called
a {\em Morita context} and the pair $(X, Y)$
a {\em Morita pair}.

If $X$ is a sup-lattice, we denote by ${\cal Q}(X)$ the
quantale of sup-preserving operators on $X$.


\section {Morita pairs between quantales}

In this section we characterize Morita pairs
between m-regular quantales.

\begin{theorem}\label{th_mor_pair}\label{theorem1}
Let $X$ and $Y$ be sup-lattices. Then $(X,Y)$ is
a Morita equivalence pair between some
m-regular quantales
if and only if there exist surjective sup-preserving maps
$p: X\otimes Y\otimes X\rightarrow X$ and
$q: Y\otimes X\otimes Y\rightarrow Y$
such that

\begin{enumerate}
\item $p(p(x_1\otimes y_1\otimes x_2)\otimes y_2\otimes x_3)=$
$p(x_1\otimes q(y_1\otimes x_2\otimes y_2)\otimes x_3)=
p(x_1\otimes y_1\otimes p(x_2\otimes y_2\otimes x_3))$;

\item $q(q(y_1\otimes x_1\otimes y_2)\otimes x_2\otimes y_3)=$
$q(y_1\otimes p(x_1\otimes y_2\otimes x_2)\otimes y_3)=$%
$q(y_1\otimes x_1\otimes q(y_2\otimes x_2\otimes y_3))$;
\item $p( - \otimes x_1) = p( - \otimes x_2)$ implies
$x_1=x_2$;
\item $p(x_1\otimes -) = p(x_2\otimes -)$ implies
$x_1=x_2$;
\item $q( - \otimes y_1) = q( - \otimes y_2)$ implies
$y_1=y_2$ and
\item $q(y_1\otimes -) = q(y_2\otimes -)$ implies
$y_1=y_2$
\end{enumerate}

\noindent for each $x_i\in X$, $y_i\in Y$ ($i=1,2,3$).
\end{theorem}

\begin{proof}
Assume that $X$ and $Y$ are sup-lattices,
$p: X\otimes Y\otimes X\rightarrow X$ and
$q: Y\otimes X\otimes Y\rightarrow Y$ are surjective sup-preserving
maps which satisfy conditions 1-6.
\vskip0.2cm
We recall that
$L_a : X\rightarrow X$ and
$R_b : Y\rightarrow Y$ are the sup-preserving operators given by
$L_a(x) = p(a\otimes x)$ and $R_b(y) = q(b\otimes y)$,
where $a\in X\otimes Y$ and
$b\in Y\otimes X$.
\vskip0.2cm

So we can think of $L_a(x)$ ($R_b(y)$)
as a left module action by $L_a$ ($R_b$) on $X$ ($Y$).
Identities 3) and 5)
state that these actions are left separating.

Recall that $L_a\in {\cal Q}(X)$
and $R_b\in {\cal Q}(Y)$.
Let $\ccl = \{L_a : a\in X\otimes Y\}$ and
$\ccr = \{R_b : b\in Y\otimes X\}$. Evidently,
$\ccl$ is a subsup-lattice of ${\cal Q}(X)$ and
$\ccr = \{R_b : b\in Y\otimes X\}$ is
 a subsup-lattice of {${\cal Q}(Y)$}.

We will only point out the \lq\lq left''
considerations in constructing
the Morita  context whenever the \lq\lq right'' ones follow
either by symmetry or in a similar way.


First note that $\ccl$ is a subquantale of ${\cal Q}(X)$.
Indeed, we have that
\begin{eqnarray*}
L_{x_1\otimes y_1}L_{x_1\otimes y_1}(z) & = &
p(x_1\otimes y_1\otimes p(x_2\otimes y_2\otimes z))\\
& = & p(p(x_1\otimes y_1\otimes x_2)\otimes y_2\otimes z)\\
& = &L_{p(x_1\otimes y_1\otimes x_2)\otimes y_2} (z)
\end{eqnarray*}
for each $z\in X$ and it
follows that $L_aL_b\in \ccl$ whenever $a,b\in X\otimes Y$.


Next observe that $\ccl$ is m-regular according to 3) and the surjectivity
of $p$. Indeed, if $x\in X$, $y\in Y$, then there are $x_i, u_i\in X$,
$y_i\in Y$, $i\in I$ such that $x=\bigvee %
p(x_i\otimes y_i\otimes u_i)$. Then we have that
\begin{eqnarray*}
\bigvee_{i\in I}L_{x_i\otimes y_i}L_{u_i\otimes y} & = &
L_{\bigvee_{i\in I} p(x_i\otimes y_i\otimes u_i)\otimes y} =
L_{x\otimes y}.
\end{eqnarray*}

So we have that $\ccl$ is an essential module over $\ccl$.

Now, let $a, b\in X\otimes Y$ and assume that
$L_a L_{x\otimes y}=L_b L_{x\otimes y}$ for all $x\in X$ and all
$y\in Y$. Then
$L_a(L_{x\otimes y}(z))=L_b(L_{x\otimes y}(z))$
for all $x, z\in X$ and all  $y\in Y$ i.e.
$L_a(p({x\otimes y}\otimes z))=L_b(p({x\otimes y}\otimes z))$
for all $x, z\in X$ and all  $y\in Y$ i.e.
$L_a(u)=L_b(u)$
for all $u\in X$ because  $p$ is a surjective map.

Similarly, let us assume that
$L_c L_a=L_c L_b$ for all $c\in X\otimes Y$. Then,
for all $z\in X$, we have that
$p(c\otimes L_a(z))=p(c\otimes L_b(z))$ i.e.
$L_a(z)=L_b(z)$ by the condition 3).

Alltogether, $\ccl$ is an m-regular quantale.


Let us endow the sup-lattice $X$ with the structure
of an $\ccl, \ccr$-bimodule, by letting
$L_a\cdot x = p(a\otimes x)$ and $x\cdot R_b = p(x\otimes b)$.
In a similar way, we will endow the sup-lattice $Y$ with the
structure of an $\ccr, \ccl$-bimodule, by letting
$R_b\cdot y = q(b\otimes y)$ and $y\cdot L_a = q(y\otimes a)$.


Define pairings $( -, -) : X\times Y \rightarrow \ccl$ and
$[ - , -] : Y\times X \rightarrow \ccr$ by setting
$(x,y) = L_{x\otimes y}$ and $[y,x] = R_{y\otimes x}$,
where $x\in X$ and $y\in Y$.

\vskip0.1cm

The fact that the pairings are balanced bimodule maps and the \lq\lq linking'' properties
between the two pairings are verified easily.

\vskip0.1cm


For the converse direction, suppose that
$(A,B,X,Y,( -, -),[ -, -])$ is a Morita context. Define
$\tilde{p}:X\times Y\times X\rightarrow X$ and
$\tilde{q}: Y\times X\times Y\rightarrow Y$ by letting
$$\tilde{p}(x_1, y, x_2) = (x_1,y)\cdot x_2$$
and
$$\tilde{q}(y_1, x, y_2) = [y_1,x]\cdot y_2.$$
It is obvious that $\tilde{p}$ and $\tilde{q}$ are
tri-sup-preserving maps.


If
$p:X\otimes Y\otimes X\rightarrow X$ and
$q: Y\otimes X\otimes Y\rightarrow Y$ are their
extensions through the sup-lattice tensor
product, then it is easily seen from the definition
of a Morita context that
Conditions 1-6 hold (since they are valid on generators)
and both $p$ and $q$
are surjective sup-preserving maps (since the maps $(-, -)$ and $[-, -]$ are
surjections onto $A$ and $B$ and both $X, Y$ are m-regular).
\end{proof}


\section{The involutive case}

In this section we consider the consequences
of Theorem 1  
when the quantales in the Morita context are involutive. If $X$
is a sup-lattice, then
we are able to form its conjugate sup-lattice $X^*$:
as a sup-lattice, we have $X^{*}=X$ and ${}^{*}:X\to X^{*}$,
$x\mapsto x^{*}$ denotes the identity map.

To each given $A, B$-bimodule $X$, $A, B$ involutive quantales, there naturally
corresponds a $B,A$-bimodule $X^{*}$ with a left
$B$-action and a right $A$ action on $X^{*}$ given by
$$
b\cdot x^{*}=(x\cdot b^{*})^{*}, x^{*}\cdot a=(a^{*}\cdot x)^{*}
$$
\noindent for all $a\in A$, $b\in B$.

Recall that an imprimitivity bimodule between two m-regular
involutive quantales $A$ and $B$ (see \cite{paseka})
is a sup-lattice $X$ which is an m-regular $A,B$-bimodule and
which is equipped with inner products
$_A\langle -, -\rangle : X\times X \rightarrow A$ and
$\langle -, -\rangle_B: X\times X\rightarrow B$,
such that $(A,X,_A\langle\cdot,\cdot\rangle)$
and $(B,X,\langle\cdot,\cdot\rangle_B)$ are
full Hilbert modules and
$_A\langle x,y\rangle z = x\langle y,z\rangle_B$ - the
compatibility condition.


If $A$ and $B$ are m-regular involutive quantales and
$X$ is an imprimitivity bimodule between
$A$ and $B$, then we are able to form an imprimitivity
$B,A$-bimodule $X^*$ in the natural way,
letting $_B\langle x^*,y^*\rangle = \langle y,x\rangle_B^*$ and
$\langle x^*,y^*\rangle_A = _A\langle y,x\rangle^*$.

Then, two m-regular involutive quantales
$A$ and $B$ are Morita equivalent via a Morita context
$(A,B,X,X^{*},(- , -),[ -, - ])$ if and only if
$X$ is an imprimitivity bimodule between
$A$ and $B$. Namely, given an imprimitivity $A, B$-bimodule, letting
$(x, y^{*})=_A\langle x, y\rangle$, $[x^{*}, y]=\langle x, y\rangle_B$,
the condition 2) from the definition of Morita equivalence
holds by the compatibility condition. Since
$_A\langle -, -\rangle$ and $\langle -, -\rangle_B$
are full, it follows that that $( -, -)$ and $[ -, -]$ are surjective.
The converse direction goes the same way.



\begin{theorem}\label{theorem2}
Let $X$ be an sup-lattice. Then $X$ is
an imprimitivity bimodule between certain
Morita equivalent m-regular involutive quantales if
and only if there exists a
surjective sup-preserving map
$p : X\otimes X^*\otimes X\rightarrow X$ such that

a) $p(p(x_1\otimes x_2^*\otimes x_3)\otimes x_4^*\otimes x_5) =
p(x_1\otimes p(x_4\otimes x_3^*\otimes x_2)^*\otimes x_5) =
p(x_1\otimes x_2^*\otimes p(x_3\otimes x_4^*\otimes x_5))$,

b)  $p( - \otimes x_1) = p( - \otimes x_2)$ implies
$x_1=x_2$ and

c)  $p( x_1 \otimes -) = p( x_2 \otimes -)$ implies
$x_1=x_2$

\noindent for each $x_i\in X$ ($i=1,2,3,4,5$).

\end{theorem}
\begin{proof}
Suppose that a sup-lattice $X$ and a map $p$ are
given such that a), b) and c) are
satisfied. Let $q:X^*\otimes X\otimes X^*\rightarrow X^*$
be the map given by
$$q(x^*\otimes y\otimes z^*) = p(z\otimes y^*\otimes x)^*,
\ \ x,y,z\in X.$$
An immediate verification shows that
$p$ and $q$ satisfy conditions 1-6 of the Theorem
\ref{theorem1}.

Note that $\ccl$ and $\ccr$ are m-regular involutive quantales with
respect to the  involutions on $\ccl$ and $\ccr$ defined by:
if $x,y\in X$, we let $L_{x\otimes y^*}^* = L_{y\otimes x^*}$ and
$R_{y^*\otimes x}^* = R_{x^*\otimes y}$.
The fact that the mappings $L\rightarrow L^*$ and $R\rightarrow R^*$ are
indeed involutions on $\ccl$ and $\ccr$ follows easily from condition a).

Thus we have shown that
$(\ccl,\ccr,X,X^{*},(- , -),[ -, - ])$ is a Morita context.

The converse direction is obtained
similarly as the converse direction in Theorem 1.%
\end{proof}

\section {Morita pairs between $\kappa$-quantales}

The following definitions are either taken  from \cite{nkuimi}
or they are a straightforward reformulation of those in section 1.
In what follows we shall always assume that $\kappa$ will be
an infinite regular cardinal.

\medskip

By a {\em $\kappa$-join semilattice} we understand a
poset $(S, \leq )$ for which every subset of cardinality
strictly less than $\kappa$ has a join. Such a join is called
a {\em $\kappa$-join}. {\em Morphisms} between $\kappa$-semilattices
are required to preserve $\kappa$-joins. Note that any
$\kappa$-join semilattice has arbitrary finite joins, especially
it has the bottom element $0$. If $\kappa=\omega$ we call
$\omega$-join semilattices only {\em join semilattices}.

By a {\em $\kappa$-quantale} is meant a $\kappa$-join semilattice
equipped with an associative multiplication; this multiplication
distributes over arbitrary $\kappa$-joins. A $\kappa$-quantale
is said to be {\em unital} if the multiplication has a unit.
If $\kappa=\omega$ we call
$\omega$-quantales {\em m-semilattices} (see \cite{pasekaz}).

An {\em involution}
on a $\kappa$-join semilattice $S$ is a unary operation  such that
$$
\begin{array}{r c l}
a^{**}&=&a,\\
(\bigvee_{i\in\kappa} a_{i})^{*}&=&\bigvee_{i\in\kappa} a_{i}^{*}
\end{array}
$$
\noindent for all $a, a_{i}\in S$, $i\in \kappa$.
An {\em involution}
on a $\kappa$-quantale $A$ is an involution on the
$\kappa$-join semilattice
$A$  such that
$$
\begin{array}{r c l}
(a\cdot b)^{*}&=&b^{*}\cdot a^{*},\\
\end{array}
$$
\noindent for all $a,  b\in A$.
A $\kappa$-join semilattice ($\kappa$-quantale) with the involution is said to be
{\em involutive}.

By a {\em morphism of} ({\em involutive}) {\em $\kappa$-quantales}
will be meant a
$\kappa$-joins- ($^{*}$-) and $\cdot$\,{}-preserving mapping
$f:A\to A'.$ If a morphism preserves the unital element
we say that it is {\em unital}.

Let $A$ be a $\kappa$-quantale.
A {\em right $\kappa$-module over} $A$ (shortly a right $A$-$\kappa$-module) is
a $\kappa$-join semilattice $M$, together with a {\em module action}
$\hbox{\rm{\_}}{\cdot}\hbox{\rm{\_}}:M\times A\to M$ satisfying

$$
\begin{array}{c l c l c}
\phantom{xxxxxxxxxx}&%
m\cdot(a\cdot b)&=&(m\cdot a)\cdot b&
\phantom{xxxxxxxxxx}\hbox{\rm (M1)}\\
\phantom{xxxxxxxxx}&%
(\bigvee X)\cdot a&=&\bigvee \{x\cdot a: x\in X\}&
\phantom{xxxxxxxxxx}\hbox{\rm ($\kappa$-M2)}\\
\phantom{xxxxxxxx}&
m\cdot \bigvee S&=&\bigvee \{m\cdot s: s\in S\}&
\phantom{xxxxxxxxxx}\hbox{\rm ($\kappa$-M3)}\\
\end{array}
$$

\noindent for all $a, b\in A$, $m\in M$, $S\subseteq A$, $|S|<\kappa$,
$X\subseteq M$, $|X|<\kappa$. Morphism of $\kappa$-modules over a
$\kappa$-quantale $A$ are morphisms of $\kappa$-join semilattices which
preserve the action by $A$.

By the same procedure as in section 1 we may define the {$\kappa$-essential}
part of a modules, {right-separatedness},  {m-regularity} and other basic
notions from the quantale case. So we shall omit it.

Two m-regular $\kappa$-quantales $A$ and $B$ are said to be {\em Morita equivalent}
if
there exist $\kappa$-join semilattices $X$ and $Y$, such that

1)$X$ is an m-regular $A,B$-$\kappa$-bimodule, $Y$ is an m-regular $B,A$-$\kappa$-bimodule;

2) there are $\kappa$-bimodule maps
$( -, -) : X\times Y\longrightarrow A$
and  $[ -, -] : Y\times X\longrightarrow B$ called pairings
such that $(x\cdot b,y) = (x, b\cdot y)$,
$[y\cdot a,x] = [y, a\cdot x]$ (that is, these maps are
{\em balanced}), $(x_1,y)\cdot x_2 = x_1\cdot [y,x_2]$,
$[y_1,x]\cdot y_2 = y_1\cdot (x,y_2)$ for each
$x,x_1,x_2\in X$, $y,y_1,y_2\in Y$, $a\in A$, $b\in B$ and

3) the $\kappa$-join-preserving maps on the  $\kappa$-join semilattice tensor
products   $X\otimes_{\kappa} Y, Y\otimes_{\kappa} X$
induced by the pairings are surjective.


The 6-tuple $(A,B,X,Y,( -, -),[ -, -])$ is called
a {\em Morita context} and the pair $(X, Y)$
a {\em Morita pair}.

If $X$ is a $\kappa$-join semilattice, we denote by ${\cal Q}_{\kappa}(X)$ the
$\kappa$-quantale of $\kappa$-join-preserving operators on $X$.

Since all arguments in the proofs of the following propositions are the same
as in Theorem \ref{theorem1} and Theorem \ref{theorem2} we shall omit them.

So we only state the corresponding theorems in the $\kappa$-setting.

\begin{theorem}\label{theorem3}
Let $X$ and $Y$ be $\kappa$-join semilattices. Then $(X,Y)$ is
a Morita equivalence pair between some
m-regular $\kappa$-quantales
if and only if there exist surjective $\kappa$-join semilattice morphisms
$p: X\otimes Y\otimes X\rightarrow X$ and
$q: Y\otimes X\otimes Y\rightarrow Y$
satisfying Conditions 1-6 from the Theorem \ref{theorem1}.
\end{theorem}

\begin{corollary}\label{corollary1}
Let $X$ and $Y$ be join semilattices. Then $(X,Y)$ is
a Morita equivalence pair between some
m-regular m-semilattices
if and only if there exist surjective join semilattice morphisms
$p: X\otimes Y\otimes X\rightarrow X$ and
$q: Y\otimes X\otimes Y\rightarrow Y$
satisfying Conditions 1-6 from the Theorem \ref{theorem1}.
\end{corollary}

\begin{theorem}\label{theorem4}
Let $X$ be a $\kappa$-join semilattice. Then $(X,X^{*})$ is
a Morita equivalence pair between some
involutive m-regular $\kappa$-quantales
if and only if there exists a surjective $\kappa$-join semilattice morphism
$p: X\otimes X^{*}\otimes X\rightarrow X$
satisfying Conditions a)-c) from the Theorem \ref{theorem2}.
\end{theorem}

\begin{corollary}\label{corollary2}
Let $X$ be a join semilattice. Then $(X,X^{*})$ is
a Morita equivalence pair between some
involutive m-regular m-semilattices
if and only if there exists a surjective join semilattice morphism
$p: X\otimes X^{*}\otimes X\rightarrow X$
satisfying Conditions a)-c) from the Theorem \ref{theorem2}.
\end{corollary}

\vskip15mm
{\parindent0pt
\begin{tabular}{@{}l l}
Author's address:  &Department of  Mathematics\\
                   &Masaryk University Brno\\
                   &Jan\'a\v{c}kovo n\'am. 2a\\
                   &662 95 Brno\\
                   &Czech Republic\\
& \\
e-mail:            &paseka@math.muni.cz
\end{tabular}
}

\end{document}